\documentclass{snmp2001}

\begin{document}

\newcommand{\HH}{{\mathcal H}}
\newcommand{\LL}{{\mathcal L}}
\newcommand{\Sc}{{\mathcal S}}

\newcommand{\RR}{{\mathbb R}}
\newcommand{\ZZ}{{\mathbb Z}}
\newcommand{\PP}{{\mathbb P}}

\newcommand{\im}{\mathop{\rm Im}\nolimits}

\newcommand{\Der}{\mathop{\rm Der}\nolimits}
\newcommand{\vth}{\vartheta}
\newcommand{\be}{\beta}
\newcommand{\na}{\nabla}
\newcommand{\om}{\omega}
\newcommand{\Ga}{\Gamma}
\newcommand{\Om}{\Omega}
\newcommand{\bom}{\bar{\omega}}
\newcommand{\bGa}{\bar{\Gamma}}
\newcommand{\bOm}{\bar{\Omega}}
\newcommand{\Z}{Z^2/\{0\}}
\newcommand{\pa}{\partial}
\newcommand{\tF}{\tilde{F}}
\newcommand{\tH}{\tilde{H}}
\newcommand{\vom}{\vec{\om}}
\newcommand{\tom}{\tilde{\omega}}
\newcommand{\tq}{\tilde{q}}
\newcommand{\tz}{\tilde{z}}
\newcommand{\tx}{\tilde{x}}
\newcommand{\ty}{\tilde{y}}
\newcommand{\ta}{\tilde{a}}
\newcommand{\tN}{\tilde{N}}
\newcommand{\tw}{\tilde{w}}
\newcommand{\J}{{\cal{J}}}
\newcommand{\F}{{\cal{F}}}
\newcommand{\R}{{\cal{R}}}
\newcommand{\cf}{\mbox{f}}
\newcommand{\E}{{\cal{E}}}
\newcommand{\N}{{\cal{N}}}
\newcommand{\cS}{{\cal{S}}}
\newcommand{\g}{\vec{g}}
\newcommand{\q}{\vec{q}}
\newcommand{\p}{\vec{p}}
\newcommand{\e}{\epsilon}
\newcommand{\te}{\tilde{\epsilon}}
\newcommand{\U}{{\cal U}}
\newcommand{\cq}{\hat{q}}
\newcommand{\A}{{\cal A}}
\newcommand{\tA}{\tilde{{\cal A}}}
\newcommand{\cO}{{\cal O}}
\renewcommand{\k}{\kappa}
\newcommand{\ga}{\gamma}
\newcommand{\ve}{\vec{e}}
\newcommand{\vv}{\vec{v}}
\newcommand{\hcS}{\hat{\cS}}
\newcommand{\tcS}{\tilde{\cS}}
\newcommand{\hS}{\hat{S}}
\newcommand{\tS}{\tilde{S}}
\newcommand{\dl}{\delta}
\newcommand{\Dl}{\Delta}
\renewcommand{\th}{\theta}
\newcommand{\ra}{\rightarrow}
\newcommand{\al}{\alpha}
\newcommand{\sg}{\sigma}
\newcommand{\Sg}{\Sigma}
\newcommand{\bM}{\bar{M}}
\newcommand{\z}{\zeta}
\newcommand{\hQ}{\hat{Q}}
\newcommand{\hv}{\hat{v}}
\newcommand{\La}{\Lambda}
\newcommand{\la}{\lambda}
\newcommand{\tla}{\tilde{\lambda}}
\newcommand{\bq}{\bar{q}}
\newcommand{\bp}{\bar{p}}
\newcommand{\bQ}{\bar{Q}}
\newcommand{\bE}{\bar{E}}
\newcommand{\rf}{\bar{f}}
\newcommand{\nid}{\noindent}
\newcommand{\rc}{S_\omega}
\newcommand{\hrc}{\hat{S}_\omega}
\newcommand{\bW}{\bar{W}}
\newcommand{\hN}{\hat{N}}
\newcommand{\hF}{\hat{F}}
\newcommand{\hk}{\hat{k}}
\newcommand{\M}{{\cal M}}
\newcommand{\D}{{\cal D}}
\newcommand{\W}{{\cal W}}
\newcommand{\C}{{\cal C}}
\newcommand{\lag}{\langle}
\newcommand{\rag}{\rangle}
\newcommand{\es}{\mbox{\boldmath$s$}}
\newcommand{\eom}{\mbox{\boldmath$\omega$}}

%\FirstPageHeading{Li}
% The parameter is the label of the article. Good choice is the last name of the first author

\ShortArticleName{Integrable Structures for 2D Euler Equations} % maximum 75 symbols

\ArticleName{Integrable Structures for 2D Euler Equations \\
of Incompressible Inviscid Fluids}

% Names of the authors for the title of the paper
\Author{Yanguang (Charles)  LI~$^\dag$}
\AuthorNameForHeading{Y. Li}
\AuthorNameForContents{LI Y.}

% Address of First Author
\Address{$^\dag$~Department of Mathematics, University of Missouri, 
Columbia, MO 65211, USA}
\EmailD{cli@math.missouri.edu}

% In the case of the same organization, please use the following standard
%\Author{First Names LASTNAME and Second COAUTHOR}
%\AuthorNameForHeading{F.N. Lastname and S. Coauthor}
%\AuthorNameForContents{LASTNAME F.N. and COAUTHOR S.}
%\Address{Address of Author(s), Country}
%\Email{email1@address, email2@address}

\Abstract{In this article, I will report a Lax pair structure, a 
B\"acklund-Darboux 
transformation, and the investigation of homoclinic 
structures for 2D Euler equations of incompressible inviscid fluids.}

\section{Introduction}

The governing equation of turbulence, that we are interested in, is the 
incompressible 2D Navier-Stokes equation under periodic boundary conditions. 
We are particularly 
interested in investigating the dynamics of 2D Navier-Stokes equation in 
the infinite Reynolds number limit and of 2D Euler equation. Our approach 
is different from many other studies on 2D Navier-Stokes equation in which 
one starts with Stokes equation to prove results on 2D Navier-Stokes 
equation for small Reynolds number. In our studies, we start with 2D Euler 
equation and view 2D Navier-Stokes equation for large Reynolds number as 
a (singular) perturbation of 2D Euler equation. 2D Euler equation is a 
Hamiltonian system with infinitely many Casimirs. To understand the nature 
of turbulence, we start with investigating the hyperbolic structure of 2D 
Euler equation. We are especially interested in investigating the 
possible homoclinic structures. 

In \cite{Li00}, we studied a linearized 2D Euler equation at a fixed point. 
The linear system decouples into infinitely many one-dimensional invariant 
subsystems. The essential spectrum of each invariant subsystem is a band 
of continuous spectrum on the imaginary axis. Only finitely many of these 
invariant subsystems have point spectra. The point spectra can be computed 
through continued fractions. Examples show that there are indeed 
eigenvalues with positive and negative real parts. Thus, there is linear 
hyperbolicity. 

In \cite{Li01a} and \cite{LY01}, a Lax pair and a B\"acklund-Darboux 
transformation were found for the 2D Euler equation. Typically, 
B\"acklund-Darboux transformation can be used to generate homoclinic 
orbits \cite{Li00a}.

The 2D Euler equation can be written in the vorticity form,
\begin{equation}
\pa_t \Om + \{ \Psi, \Om \} = 0 \ ,
\label{li:euler}
\end{equation}
where the bracket $\{\ ,\ \}$ is defined as
\[
\{ f, g\} = (\pa_x f) (\pa_y g) - (\pa_y f) (\pa_x g) \ ,
\]
where $\Psi$ is the stream function given by,
\[
u=- \pa_y \Psi \ ,\ \ \ v=\pa_x \Psi \ ,
\]
$u$ and $v$ are the velocity components, and the relation between 
vorticity $\Om$ and stream 
function $\Psi$ is,
\[
\Om =\pa_x v - \pa_y u =\Dl \Psi \ .
\]

\section{A Lax Pair and a Darboux Transformation}

\begin{theorem}[Li, \cite{Li01a}]
The Lax pair of the 2D Euler equation (\ref{li:euler}) is given as
\begin{equation}
\left \{ \begin{array}{l} 
L \varphi = \la \varphi \ ,
\\
\pa_t \varphi + A \varphi = 0 \ ,
\end{array} \right.
\label{li:laxpair}
\end{equation}
where
\[
L \varphi = \{ \Om, \varphi \}\ , \ \ \ A \varphi = \{ \Psi, \varphi \}\ ,
\]
and $\la$ is a complex  constant, and $\varphi$ is a complex-valued function.
\label{li:2dlp}
\end{theorem}
Consider the Lax pair (\ref{li:laxpair}) at $\la =0$, i.e.
\begin{eqnarray}
& & \{ \Om, p \} = 0 \ , \label{li:d1} \\
& & \pa_t p + \{ \Psi, p \} = 0 \ , \label{li:d2} 
\end{eqnarray}
where we replaced the notation $\varphi$ by $p$.
\begin{theorem}[Li and Yurov, \cite{LY01}]
Let $f = f(t,x,y)$ be any fixed solution to the system 
(\ref{li:d1}, \ref{li:d2}), we define the Gauge transform $G_f$:
\begin{equation}
\tilde{p} = G_f p = \frac {1}{\Om_x}[p_x -(\pa_x \ln f)p]\ ,
\label{li:gauge}
\end{equation}
and the transforms of the potentials $\Om$ and $\Psi$:
\begin{equation}
\tilde{\Psi} = \Psi + F\ , \ \ \ \tilde{\Om} = \Om + \Dl F \ ,
\label{li:ptl}
\end{equation}
where $F$ is subject to the constraints
\begin{equation}
\{ \Om, \Dl F \} = 0 \ , \ \ \ \{ \Om +\Dl F, F \} = 0\ .
\label{li:constraint}
\end{equation}
Then $\tilde{p}$ solves the system (\ref{li:d1}, \ref{li:d2}) at 
$(\tilde{\Om}, \tilde{\Psi})$. Thus (\ref{li:gauge}) and 
(\ref{li:ptl}) form the Darboux transformation for the 2D 
Euler equation (\ref{li:euler}) and its Lax pair (\ref{li:d1}, \ref{li:d2}).
\label{li:dt}
\end{theorem}

\section{Preliminaries on Linearized 2D Euler Equation}

We consider the two-dimensional incompressible Euler equation
written in vorticity form(\ref{li:euler})
under periodic boundary conditions in both $x$ and $y$ directions
with period $2\pi$. We also require that both 
$u$ and $v$ have means zero,
\[
\int_0^{2\pi}\int_0^{2\pi} u\ dxdy =\int_0^{2\pi}\int_0^{2\pi} v\ dxdy=0.
\]
\nid
We expand $\Om$ into Fourier series,
\[
\Om =\sum_{k\in Z^2/\{0\}} \om_k \ e^{ik\cdot X}\ ,
\]
\nid
where $\om_{-k}=\overline{\om_k}\ $, $k=(k_1,k_2)^T$, 
$X=(x,y)^T$. In this paper, we confuse $0$ with $(0,0)^T$, the context 
will always make it clear. By the relation between vorticity $\Om$ and stream 
function $\Psi$,
the system (\ref{li:euler}) can be rewritten as the following kinetic system,
\begin{equation}
\dot{\om}_k = \sum_{k=p+q} A(p,q) \ \om_p \om_q \ ,
\label{li:Keuler}
\end{equation}
\nid
where $A(p,q)$ is given by,
\begin{equation}
A(p,q)= \frac{1}{2}[|q|^{-2}-|p|^{-2}](p_1 q_2 -p_2 q_1)\ ,
\label{li:Af} 
\end{equation}
\nid
where $|q|^2 =q_1^2 +q_2^2$ for $q=(q_1,q_2)^T$, similarly for $p$.

\begin{figure}[ht]
  \begin{center}
    \leavevmode
      \setlength{\unitlength}{2ex}
  \begin{picture}(36,27.8)(-18,-12)
    \thicklines
\put(0,-14){\vector(0,1){28}}
\put(-18,0){\vector(1,0){36}}
\put(0,15){\makebox(0,0){$\Im \{ \la \}$}}
\put(18.5,0){\makebox(0,0)[l]{$\Re \{ \la \}$}}
\put(0.1,-7){\line(0,1){14}}
\put(.2,-.2){\makebox(0,0)[tl]{$0$}}
\put(-0.2,-7){\line(1,0){0.4}}
\put(-0.2,7){\line(1,0){0.4}}
\put(2.0,-6.4){\makebox(0,0)[t]{$-i2|b|$}}
\put(2.0,7.6){\makebox(0,0)[t]{$i2|b|$}}
\put(2.4,3.5){\circle*{0.5}}
\put(-2.4,3.5){\circle*{0.5}}
\put(2.4,-3.5){\circle*{0.5}}
\put(-2.4,-3.5){\circle*{0.5}}
\put(5,4){\circle*{0.5}}
\put(-5,4){\circle*{0.5}}
\put(5,-4){\circle*{0.5}}
\put(-5,-4){\circle*{0.5}}
\put(8,6){\circle*{0.5}}
\put(-8,6){\circle*{0.5}}
\put(8,-6){\circle*{0.5}}
\put(-8,-6){\circle*{0.5}}
\end{picture}
  \end{center}
\caption{The spectrum of $\LL_{\hk}$.}
\label{li:spla2}
\end{figure}
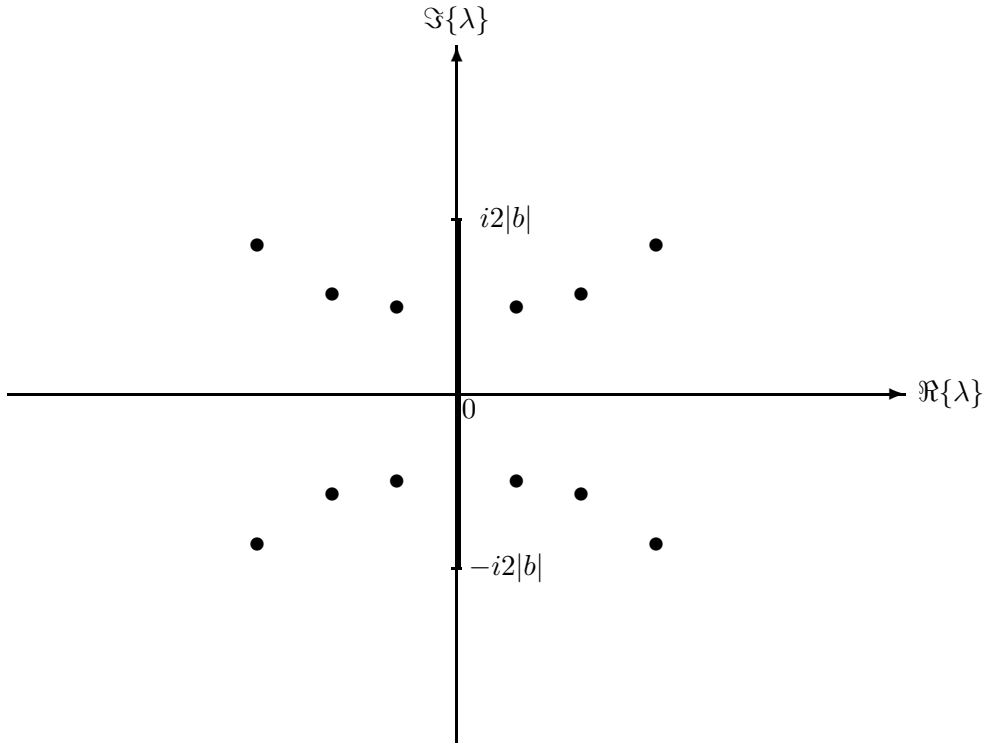

We denote $\{ \om_k \}_{k\in \Z}$ by $\om$. For any fixed 
$p \in Z^2/\{0\}$, we consider the simple fixed point 
$\om^*$:
\begin{equation}
\om^*_p = \Ga,\ \ \ \om^*_k = 0 ,\ \mbox{if} \ k \neq p \ \mbox{or}\ -p,
\label{li:fixpt}
\end{equation}
\nid
of the 2D Euler equation (\ref{li:Keuler}), where 
$\Ga$ is an arbitrary complex constant. 
The {\em{linearized two-dimensional Euler equation}} at $\om^*$ is given by,
\begin{equation}
\dot{\om}_k = A(p,k-p)\ \Ga \ \om_{k-p} + A(-p,k+p)\ \bar{\Ga}\ \om_{k+p}\ .
\label{li:LE}
\end{equation}
\begin{definition}[Classes]
For any $\hk \in \Z$, we define the class $\Sg_{\hk}$ to be the subset of 
$\Z$:
\[
\Sg_{\hk} = \bigg \{ \hk + n p \in \Z \ \bigg | \ n \in Z, \ \ p \ \mbox{is 
specified in (\ref{li:fixpt})} \bigg \}.
\]
\label{li:classify}
\end{definition}
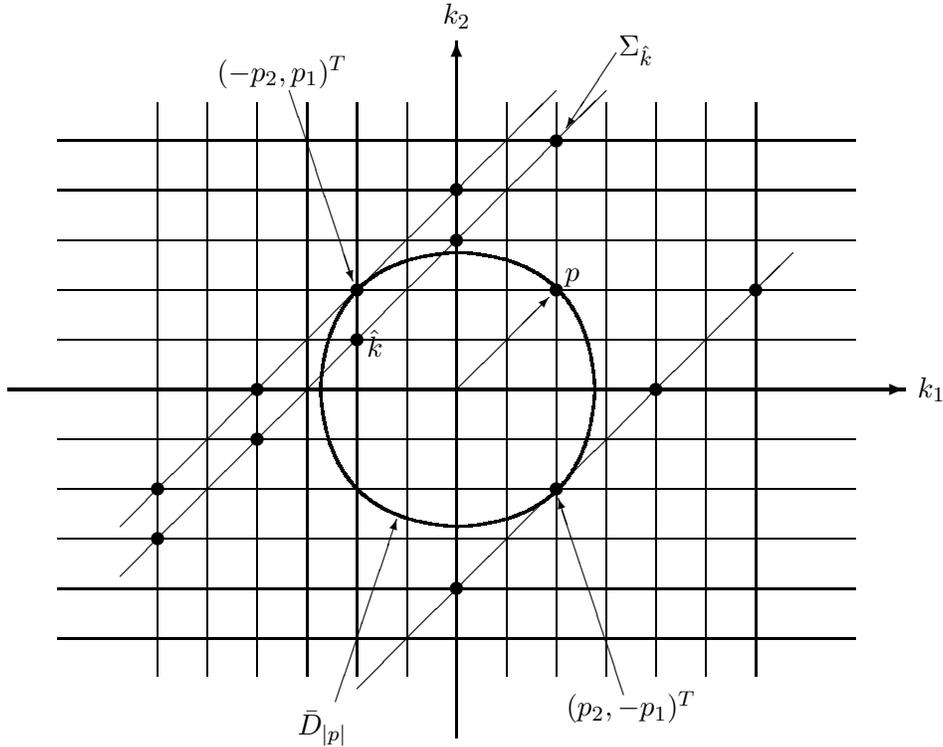
\begin{figure}[ht]
  \begin{center}
    \leavevmode
      \setlength{\unitlength}{2ex}
  \begin{picture}(36,27.8)(-18,-12)
%%  \begin{picture}(40,310)(-20,-150)
%%%    \linethickness{5pt}
    \thinlines
\multiput(-12,-11.5)(2,0){13}{\line(0,1){23}}
\multiput(-16,-10)(0,2){11}{\line(1,0){32}}
    \thicklines
\put(0,-14){\vector(0,1){28}}
\put(-18,0){\vector(1,0){36}}
\put(0,15){\makebox(0,0){$k_2$}}
\put(18.5,0){\makebox(0,0)[l]{$k_1$}}
%
%% \put(0,0){\circle{100}}
%% circle (because latex will only draw circles upto 40pts)
%
% this circle a little too small
%\qbezier(-5,0)(-4.775,4.775)(0,5)
%\qbezier(0,5)(4.775,4.775)(5,0)
%\qbezier(5,0)(4.775,-4.775)(0,-5)
%\qbezier(0,-5)(-4.775,-4.775)(-5,0)
%
%\put(3.65,3.65){\circle*{0.5}}
%\put(0,0){\vector(1,1){3.5}}
%\put(4.5,3.5){$P$}
%
%\put(3.65,-3.65){\circle*{0.5}}
%
% this circle hits required points
\qbezier(-5.5,0)(-5.275,5.275)(0,5.5)
\qbezier(0,5.5)(5.275,5.275)(5.5,0)
\qbezier(5.5,0)(5.275,-5.275)(0,-5.5)
\qbezier(0,-5.5)(-5.275,-5.275)(-5.5,0)
    \thinlines
\put(4,4){\circle*{0.5}}
\put(0,0){\vector(1,1){3.7}}
\put(4.35,4.35){$p$}
\put(4,-4){\circle*{0.5}}
\put(8,0){\circle*{0.5}}
\put(-8,0){\circle*{0.5}}
\put(-8,-2){\circle*{0.5}}
\put(-12,-4){\circle*{0.5}}
\put(-12,-6){\circle*{0.5}}
\put(-4,2){\circle*{0.5}}
\put(-4,4){\circle*{0.5}}
\put(0,6){\circle*{0.5}}
\put(0,8){\circle*{0.5}}
\put(4,10){\circle*{0.5}}
\put(12,4){\circle*{0.5}}
\put(0,-8){\circle*{0.5}}
\put(-4,-12){\line(1,1){17.5}}
\put(-13.5,-7.5){\line(1,1){19.5}}
\put(-13.5,-5.5){\line(1,1){17.5}}
\put(-3.6,1.3){$\hat{k}$}
\put(-7,12.1){\makebox(0,0)[b]{$(-p_2, p_1)^T$}}
%\put(-7.85,11.75){\vector(1,-2){3.65}}
\put(-6.7,12){\vector(1,-3){2.55}}
\put(6.5,13.6){\makebox(0,0)[l]{$\Sg_{\hat{k}}$}}
\put(6.4,13.5){\vector(-2,-3){2.0}}
\put(7,-12.1){\makebox(0,0)[t]{$(p_2, -p_1)^T$}}
\put(6.7,-12.25){\vector(-1,3){2.62}}
\put(-4.4,-13.6){\makebox(0,0)[r]{$\bar{D}_{|p|}$}}
\put(-4.85,-12.55){\vector(1,3){2.45}}
  \end{picture}
  \end{center}
\caption{An illustration of the classes $\Sg_{\hk}$ and the disk 
$\bar{D}_{|p|}$.}
\label{li:class}
\end{figure}
\nid
See Figure \ref{li:class} for an illustration. According to the classification 
defined in Definition \ref{li:classify}, the linearized two-dimensional Euler 
equation (\ref{li:LE}) decouples into infinitely many {\em{invariant subsystems}}:
\begin{eqnarray}
\dot{\omega}_{\hat{k} + np} &=& A(p, \hat{k} + (n-1) p) 
     \ \Gamma \ \omega_{\hat{k} + (n-1) p} \nonumber \\  \label{li:CLE}\\
& & + \ A(-p, \hat{k} + (n+1)p)\ 
     \bar{\Gamma} \ \omega_{\hat{k} +(n+1)p}\ . \nonumber
\end{eqnarray}
\begin{theorem}
The eigenvalues of the linear operator $\LL_{\hk}$ defined by the 
right hand side 
of (\ref{li:CLE}), are of 
four types: real pairs ($c, -c$), purely imaginary pairs ($id, -id$), 
quadruples ($\pm c \pm id$), and zero eigenvalues.
\end{theorem}
\nid
The eigenvalues can be computed through continued fractions.
\begin{definition}[The Disk]
The disk of radius $\left| p \right|$ in $Z^2 / \left\{ 0
\right\}$, denoted by
$\bar{D}_{\left| p \right|}$, is defined as
\[ 
 \bar{D}_{\left| p \right|} = \bigg \{ k \in Z^2/ \left\{ 0 \right\} \ \bigg| 
     \ \left| k \right| \leq \left| p \right| \bigg \} \, .
\]
\end{definition}
\begin{theorem}[The Spectral Theorem] We have the following claims on 
the spectra of the linear operator $\LL_{\hk}$:
\begin{enumerate}
\item If $\Sg_{\hat{k}} \cap \bar{D}_{|p|} = \emptyset$, then the entire
$\ell_2$ spectrum of the linear operator $\LL_{\hk}$ 
is its continuous spectrum. See Figure \ref{li:spla2}, where
$b= - \frac{1}{2}|\Gamma | |p|^{-2} 
\left|
  \begin{array}{cc}
p_1 & \hat{k}_1 \\
p_2 & \hat{k}_2
  \end{array}
\right| \ .$
\item If $\Sg_{\hat{k}} \cap \bar{D}_{|p|} \neq \emptyset$, then the entire
essential $\ell_2$ spectrum of the linear operator $\LL_{\hk}$ is its 
continuous spectrum. 
That is, the residual 
spectrum of $\LL_{\hk}$ is empty, $\sg_r (\LL_{\hk}) = \emptyset$. The point 
spectrum of $\LL_{\hk}$ is symmetric with respect to both real and 
imaginary axes. 
See Figure \ref{li:spla2}.
\end{enumerate}
\label{li:spthla}
\end{theorem}

\section{A Galerkin Truncation}

\begin{figure}[ht]
  \begin{center}
    \leavevmode
      \setlength{\unitlength}{2ex}
  \begin{picture}(36,27.8)(-18,-12)
    \thinlines
\multiput(-12,-11.5)(2,0){13}{\line(0,1){23}}
\multiput(-16,-10)(0,2){11}{\line(1,0){32}}
    \thicklines
\put(0,-14){\vector(0,1){28}}
\put(-18,0){\vector(1,0){36}}
\put(0,15){\makebox(0,0){$k_2$}}
\put(18.5,0){\makebox(0,0)[l]{$k_1$}}
\qbezier(-2.75,0)(-2.6375,2.6375)(0,2.75)
\qbezier(0,2.75)(2.6375,2.6375)(2.75,0)
\qbezier(2.75,0)(2.6375,-2.6375)(0,-2.75)
\qbezier(0,-2.75)(-2.6375,-2.6375)(-2.75,0)
    \thinlines
\put(2,2){\circle*{0.5}}
\put(0,0){\vector(1,1){1.85}}
\put(2.275,2.275){$p$}
\put(-12,-10){\circle*{0.5}}
\put(-10,-8){\circle*{0.5}}
\put(-8,-6){\circle*{0.5}}
\put(-6,-4){\circle{0.5}}
\put(-4,-2){\circle*{0.5}}
\put(-2,0){\circle*{0.5}}
\put(0,2){\circle*{0.5}}
\put(2,4){\circle*{0.5}}
\put(4,6){\circle{0.5}}
\put(6,8){\circle*{0.5}}
\put(8,10){\circle*{0.5}}
\put(-14,-12){\line(1,1){24}}
\put(-5.6,-5.4){$\hat{k}$}
\put(-4.4,-13.6){\makebox(0,0)[r]{$\bar{D}_{|p|}$}}
\put(-4.85,-12.55){\vector(1,3){3.4}}
\end{picture}
\end{center}
\caption{The collocation of the modes in the Galerkin truncation.}
\label{li:model}
\end{figure}
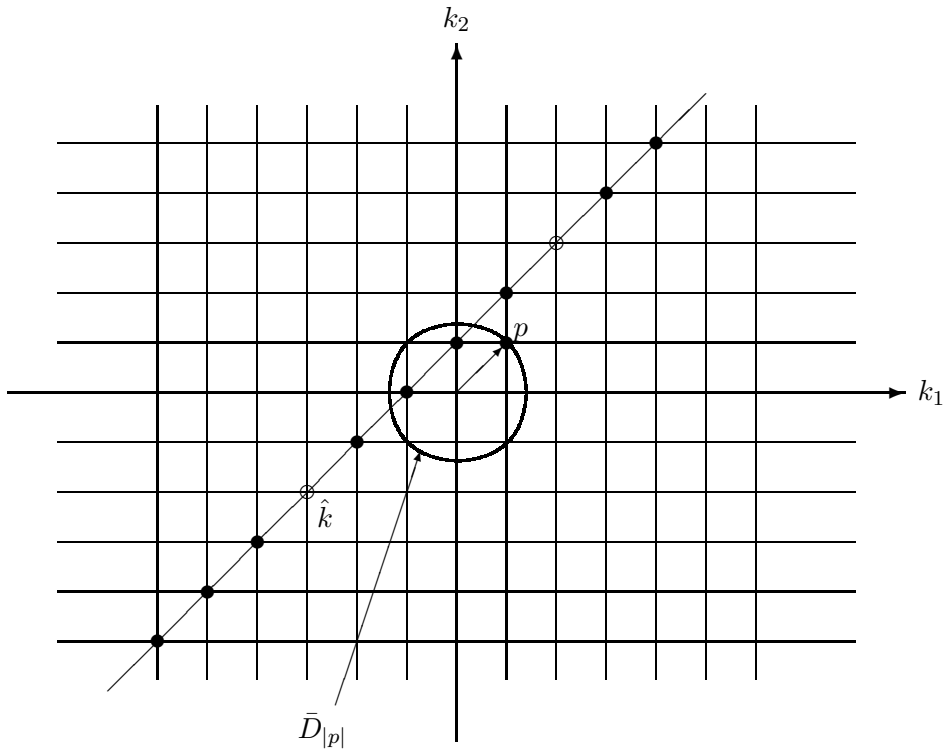

To simplify our study, we study only the case when $\om_k$ is real, $\forall 
k \in \Z$, i.e. we only study the cosine transform of the vorticity, 
\[
\Om = \sum_{k \in \Z} \om_k \cos (k \cdot X)\ ,
\]
and the 2D Euler equation (\ref{li:euler};\ref{li:Keuler}) preserves the cosine 
transform. To further simplify our study, we will study a concrete 
line of fixed points (\ref{li:fixpt}) with the 
mode $p=(1,1)^T$ parametrized by $\Ga$.
When $\Ga \neq 0$, each fixed point has $4$ eigenvalues which form a 
quadruple. These four eigenvalues appear in the only unstable invariant 
linear subsystem labeled by $\hk = (-3,-2)^T$. See Figure \ref{li:model}
for an illustration.
\nid
We computed the eigenvalues 
through continued fractions, one of them is \cite{Li00}:
\begin{equation}
\tla=2 \lambda / | \Gamma | = 0.24822302478255 \ + \ i \ 0.35172076526520\ .
\label{li:evun}
\end{equation}
We hope that a Galerkin truncation with a small 
number of modes including those inside the disk 
$\bar{D}_{\left| p \right|}$ can capture the eigenvalues.
We propose the Galerkin truncation to the linear system 
(\ref{li:CLE}) with the four modes $\hk +p$,  $\hk +2p$, $\hk +3p$, and $\hk +4p$,
\begin{eqnarray*}
\dot{\omega}_{1} &=& -A_{2} \Gamma \omega_{2} \, , \\
\dot{\omega}_{2} &=& A_{1} \Gamma \omega_{1}
- A_{3} \Gamma \omega_{3}\, , \\
\dot{\omega}_{3} &=& A_{2} \Gamma \omega_{2}
- A_{4} \Gamma \omega_{4}\, , \\
\dot{\omega}_{4} &=& A_{3} \Gamma \omega_{3} \, . 
\end{eqnarray*}
From now on, the abbreviated notations,
\begin{equation}
  \omega_n = \omega_{\hat{k}+np} \, , \ \  
  A_n = A(p,\hat{k}+np) \, , \ \  
  A_{m,n} = A(\hat{k}+mp,\hat{k}+np) \, ,  
\label{li:abbn}
\end{equation}
will be used. The eigenvalues of this four dimensional system can be 
easily calculated. It turns out that this system has a quadruple of 
eigenvalues:
\begin{equation}
\lambda = \pm \frac{\Gamma}{2 \sqrt{10}} \sqrt{1 \pm i \sqrt{35}} 
\dot{=} \pm\left( \frac{\Gamma}{2} \right) \times 0.7746 
  \times e^{\pm i \theta_1} \, , \label{li:evu}
\end{equation}
where $\theta_1 = \arctan (0.845)$, in comparison with the
quadruple of eigenvalues (\ref{li:evun}), where
\begin{displaymath}
  \lambda \dot{=} \pm \left( \frac{\Gamma}{2} \right) \times 0.43 \times 
  e^{\pm i \theta_2} \, ,
\end{displaymath}
and $\theta_2 = \arctan (1.418)$. Thus, {\em{ the quadruple of eigenvalues of
the original system is recovered by the four-mode truncation}}. We 
further study the corresponding Galerkin truncation of 2D Euler equation:
\begin{eqnarray}
\dot{\omega}_1 &=& -A_2 \ \omega_p \ \omega_2 \ , \nonumber\\
\dot{\omega}_2 &=& A_1 \ \omega_p \ \omega_1 -A_2 \ \omega_p \ \omega_3
\ , \nonumber \\
\dot{\omega}_3 &=& A_2 \ \omega_p \ \omega_2 -A_1 \ \omega_p \ \omega_4
\, , \label{li:invu} \\
\dot{\omega}_4 &=& A_2 \ \omega_p  \ \omega_3 \, , \nonumber \\
\dot{\omega}_p &=& A_{1,2} \ (\omega_3 \ \omega_4 - \omega_1
\ \omega_2) \, , \nonumber 
\end{eqnarray}
and the equations for the decoupled variables $\omega_0$ and $\omega_5$
are given by,
\begin{eqnarray*}
\dot{\omega}_0 &=& -A_1 \ \omega_p \ \omega_1 \, , \\
\dot{\omega}_5 &=& A_1 \ \omega_p \ \omega_4 \, .
\end{eqnarray*}
where
\begin{eqnarray*}
& & A_1 = - \frac{3}{10}\ , \  A_2 = \frac{1}{2} \ , \ 
  A_3 = A_2 \ , \ A_4 = A_1 \ ,  \\
& & A_{1,2} = A_1 -A_2 = - \frac{4}{5} \ , \ 
  A_{2,3} =0 \ , \ A_{3,4} = -A_{1,2} \ ;  
\end{eqnarray*}
There are three invariants for the system (\ref{li:invu}):
\begin{eqnarray}
I &=& 2 A_{1,2} (\omega_1 \omega_3 + \omega_2 \omega_4 )
        + A_2 \omega^2_p \, , \label{li:Iinv} \\[1ex]
U &=&  A_1 (\omega^2_1+ \omega^2_4 ) 
        + A_2 ( \omega^2_2 + \omega^2_3 )  \, , \label{li:Uinv} \\[1ex]
J &=&   \omega^2_p+ \omega^2_1  
        + \omega^2_2 + \omega^2_3 + \omega^2_4  \, . \label{li:Jinv}
\end{eqnarray}
$J$ is the enstrophy, and $U$ is a linear combination of the kinetic 
energy and the enstrophy. $I$ is an extra invariant which is peculiar 
to this invariant subsystem. With $I$, the explicit formula for the 
hyperbolic structure can be computed.

The common level set of these three invariants which is connected 
to the fixed point (\ref{li:fixpt}) determines the stable and
unstable manifolds of the fixed point and its negative
$-\omega^*$:
\begin{equation}
\omega_p = - \Gamma \, , \, \omega_n=0 \quad (n \in Z) \, .
\label{li:nfixpt}
\end{equation}
Using the polar coordinates:
\begin{displaymath}
  \omega_1 = r \cos \theta \, , \, 
  \omega_4 = r \sin \theta \, ; \, 
  \omega_2 = \rho \cos \vartheta \, , \, 
  \omega_3 = \rho \sin \vartheta \, ;
\end{displaymath}
we have the following explicit expressions for the stable and
unstable manifolds of the fixed point (\ref{li:fixpt}) and its
negative (\ref{li:nfixpt}) represented through the homoclinic orbits
asymptotic to the line of fixed points:
\begin{eqnarray}
\omega_p &=& \Gamma \ \tanh \tau \, , \nonumber \\
r &=& \sqrt{ \frac{A_2}{A_2-A_1}}\, \ \Gamma \ \mbox{sech}\ \tau \, , \nonumber  \\[1ex]
\theta &=& - \ \frac{A_2}{2\k} \ \mbox{ln}\ \cosh \tau + \theta_0 \, , 
\label{li:exus}\\[1ex]
\rho &=& \sqrt{\frac{-A_1}{A_2}} \  r \, ,\nonumber  \\[1ex]
\theta + \vartheta &=& \left\{
  \begin{array}{ll}
    - \arcsin \left[ \frac{1}{2} \sqrt{\frac{A_2}{-A_1}}\,  \right] 
        \ , & (\k>0) \ , \\[2ex]
\pi + \arcsin  \left[ \frac{1}{2} \sqrt{\frac{A_2}{-A_1}} \right]
\, , & (\k<0) \, ,
  \end{array} \right.  \nonumber 
\end{eqnarray}
where $A_1$ and $A_2$ are given in (\ref{li:invu}), $\tau = \k \Gamma t 
+ \tau_0$, $(\tau_0, \theta_0)$ are the two parameters
parametrizing the two-dimensional stable (unstable) manifold,
and
\begin{displaymath}
  \k = \sqrt{-A_1 A_2} \cos (\theta + \vartheta) 
   = \pm \sqrt{-A_1 A_2} \sqrt{1+ \frac{A_2}{4A_1}} \ .
\end{displaymath}
The two auxilliary variables $\om_0$ and $\om_5$ have the
expressions:
\begin{eqnarray*}
  \omega_0 &=& \frac{\alpha \beta}{1+ \beta^2}\ \mbox{sech}\ \tau \left\{ 
            \sin [ \beta \ \mbox{ln}\ \cosh \tau +  \theta_0]
            - \frac{1}{\beta} \cos
            [ \beta \ \mbox{ln}\ \cosh \tau + \theta_0 ] \right\} \,
   , \\
  \omega_5 &=& \frac{\alpha \beta}{1+ \beta^2} \ \mbox{sech}\ \tau
  \left\{ \cos [ \beta \ \mbox{ln}\ \cosh \tau + \theta_0 ]
    + \frac{1}{\beta} \sin [ \beta \ \mbox{ln}\ \cosh \tau + \theta_0]
    \right\} \, ,
\end{eqnarray*}
where
\begin{displaymath}
  \alpha = -A_1 \Gamma \k^{-1} \sqrt{\frac{A_2}{A_2-A_1}} \ , \ \ 
  \beta = - \frac{A_2}{2\k} \ .
\end{displaymath}
The graphs of these homoclinic orbits are spirals on a 2D ellipsoid,
with turning points.

\section{Conclusion}

Certain newly developed results on 2D Euler equation have been 
discussed, which include a Lax pair, a Darboux transformation, 
and the investigation on homoclinic structures.


\begin{thebibliography}{99}

\bibitem{Li00} Y. Li, \emph{On 2D Euler equations: Part I.
On the energy-Casimir stabilities and the 
spectra for linearized 2D Euler equations},
J. Math. Phys., 41, No.2, pp.~728--758, (2000).

\bibitem{Li01a} Y. Li, \emph{A Lax pair for the two dimensional 
Euler equation},
J. Math. Phys., 42, No.8, (2001).

\bibitem{LY01} Y. Li and A. Yurov, \emph{Lax pairs and Darboux transformations
for Euler equations},
Submitted, (2001).

\bibitem{Li00a} Y. Li, \emph{B\"acklund-Darboux transformations and Melnikov 
analysis for Davey-Stewartson II equations}, J. Nonlinear Sci. \textbf{10} 
2000, pp.~103--131.

\bibitem{Li01b} Y. Li, \emph{On 2D Euler equations: Part II.
Lax pairs and homoclinic structures},
Submitted, (2001).


\end{thebibliography}
\end{document}